%% file: main.tex
\documentclass{article}

\usepackage{hyperref}
\usepackage{theorem}
\usepackage{amssymb}

\input{macros}

\newcommand\pre{\mr{PreFP}}
\newcommand\post{\mr{PostFP}}
\renewcommand\fp{\mr{FP}}

\begin{document}

\title{A point on fixpoints in posets}

\author{Fr\'ed\'eric Blanqui (INRIA)}

\date{24 January 2014}

\maketitle

Let $(X,\le)$ be a {\em non-empty strictly inductive poset}, that is,
a non-empty partially ordered set such that every {\em non-empty
  chain} $Y$ has a least upper bound $\lub(Y)\in X$, a chain being a
subset of $X$ totally ordered by $\le$.

We are interested in sufficient conditions such that, given an element
$a_0\in X$ and a function $f:X\a X$, there is some ordinal $k$ such
that $a_{k+1}=a_k$, where $a_k$ is the transfinite sequence of
iterates of $f$ starting from $a_0$ (implying that $a_k$ is a fixpoint
of $f$):

\begin{itemize}\itemsep=0mm
\item $a_{k+1}=f(a_k)$
\item $a_l=\lub\{a_k\mid k<l\}$ if $l$ is a limit ordinal, \ie $l=\lub(l)$
\end{itemize}

This note summarizes known results about this problem and provides a
slight generalization of some of them.

\section{Definitions}

Let AC denote the axiom of choice\footnote{Implicitly used by Cantor
  and others, but first explicitly mentioned (and rejected) by Peano
  in \cite{peano90ma}.} \cite{zermelo04ma}.

Let $\bO$ be the class of von Neumann ordinals \cite{neumann23as}.

Let $A=\{x\in X\mid\ex k\in\bO,x=a_k\}$ be the set of all transfinite
iterates of $f$ from $a_0$.

Let $N$ be the smallest subset of $X$ containing $a_0$ and closed by
$f$ ($Z\sle X$ is closed by $f$ if, for all $z\in Z$, $f(z)\in Z$) and
non-empty $\lub$'s ($Z\sle X$ is closed by non-empty $\lub$'s if, for
all $P\sle Z$, if $P\neq\vide$, then $\lub(P)$ exists and belongs to
$Z$).

Given $Z\sle X$, let $\glb(Z)$ be, if it exists, the greatest lower
bound of $Z$.

Let $\pre(f)=\{x\in X\mid x\le f(x)\}$ be the set of {\em
  pre-fixpoints} of $f$, $\post(f)=\{x\in X\mid x\ge f(x)\}$ be the
set of {\em post-fixpoints} of $f$, and $\fp(f)=\pre(f)\cap\post(f)$
be the set of fixpoints of $f$.

$f$ is (resp. {\em strictly}) {\em monotone} (or {\em isotone} or {\em
  increasing}) (opposites: {\em anti-monotone}, {\em antitone} or {\em
  decreasing}) if $f(x)\le f(y)$ whenever $x\le y$ (resp. $f(x)<f(y)$
whenever $x<y$).

$f$ is (resp. {\em strictly}) {\em extensive} (or {\em inflationary}
or {\em progressive}) (opposites: {\em reductive}, {\em deflationary}
or {\em regressive}) if, for all $x$, $f(x)\ge x$ (resp. $f(x)>x$).

Note that: if $X$ has a least element $\bot$, then it is a
pre-fixpoint of $f$; if $f$ is extensive, then {\em every} element of
$X$ is a pre-fixpoint of $f$; if $f$ is monotone, then $\pre(f)$ and
$\post(f)$ are closed by $f$ (pre-fixpoints are mapped to
pre-fixpoints, and post-fixpoints to post-fixpoints).

\section{Searching for the origins}

We found useful references in
\cite{lassez82ipl,campbell78hm,rubin63book,rubin85book,gabbay12book,felscher62am}. The
oldest result that we could find is:\\

In 1909, Hessenberg proved that, if $X$ is a set of sets, $\le$ is
inclusion, and $f$ is extensive, then $f$ has a fixpoint
\cite{hessenberg09crelle}. This is explained in \cite{felscher62am}
(as mentioned in \cite{felscher00mail}\footnote{\cite{felscher00mail}:
  "Returning to Hessenberg, his paper

   Kettentheorie und Wohlordnung. Crelle 135 (1909) 81-133

can hardly be underestimated in its importance. Not that it was
understood by his contemporaries. But Hessenberg, analyzing Zermelos
second proof of the well ordering theorem, studied the general ways to
construct well ordered subsets of ordered sets - with the one
restriction that order always was inclusion and ordered sets were
subfamilies of power sets. In the course of this, Hessenberg stated
and proved the fixpoint theorem which thirty years later was
rediscovered - for ordered sets now - by Nicolas Bourbaki. The amazing
thing is that Hessenberg's proof is precisely the same as that given
by Bourbaki ! (only that at one small point a simpler argument can be
used due to the circumstance that Hessenberg's order is
inclusion). For details, I refer to my article in Archiv d.Math. 13
(1962) 160-165 and to my book Naive Mengen und Abstrakte Zahlen from
1979 , p.200 ff ."} and \cite{gabbay12book} p. 158). Unfortunately,
these two papers are in German and I cannot get a more precise idea
of their contents.\\

In 1922, Kuratowski proved that, if $X$ is a set of sets, $\le$ is
inclusion, and $f$ is extensive, then $f$ has a fixpoint, namely
$\lub(N)$ \cite{kuratowski22fm}. He also proved that $N$ is equal to
$A$ and that, if $f$ is also monotone, then $\lub(N)$ is the smallest
fixpoint of $f$ that is greater than or equal to $a_0$.\\

In 1927, Knaster and Tarski proved that, if $X$ is a set of sets,
$\le$ is inclusion and $f$ is monotone, then $f$ has a fixpoint
\cite{knaster28aspm}.\\

In 1939, Tarski extended this result to arbitrary complete lattices
(every subset, empty or not, has a $\lub$ and a $\glb$)
\cite{tarski55pjm}. In fact, he proved that $\fp(f)$ itself is a
non-empty complete lattice. In 1951, Davis proved the converse, that
is, a lattice is complete if every monotone map has a fixpoint
\cite{davis55pjm}. By the way, note that, in 1941, Frink proved that
any complete lattice is compact in the interval topology
\cite{frink42tams}.

On the other hand, Tarski did not study whether the least fixpoint of
$f$, proved to be $\glb(\post(f))$, can be reached by transfinite
iteration. However, he mentioned p. 305 that, if $X$ is an
$\w$-complete lattice (every countable subset has a $\lub$), $f$ is
$\w$-continuous ($f(\lub Y)=\lub f(Y)$ for every countable $Y\sle X$)
and $a_0$ is the least element $\bot$ of $X$, then $a_{\w+1}=a_\w$, a
result sometimes attributed to Kleene because it is indeed used in his
proof of the {\em first} recursion theorem with $X$ being the set of
partial functions in \cite{kleene52book} p. 348.\\

At the end of the 30's, Bourbaki\footnote{A result due to Chevalley
  after \cite{campbell78hm}.} proved that, if $X$ is a strictly
inductive poset and $f$ is extensive, then $f$ has a fixpoint
\cite{bourbaki39book,bourbaki49am}. This is a generalization of
Hessenberg and Kurakowski's results. This result was rediscovered or
proved (because Bourbaki published a proof in 1949 only) at the end of
the 40's or beginning of the 50's by many other people:
\cite{kneser50mz,szele50pm,witt50rmh,witt50mn,vaughan52pjm,inagaki52mjou}.\\

Note that fixpoint theorems assuming that $f$ is extensive easily
follow from Zorn's maximal principle (equivalent to AC) saying that
any non-empty inductive (= chain-complete) set of sets ordered by
inclusion (every $\sle$-chain, including the empty one, has a $\lub$)
has a maximal element (for inclusion) \cite{zorn35bams}, or Tukey's
maximal principle (equivalent to AC too) generalizing Zorn's one to
arbitrary inductive posets (hence saying that any non-empty inductive
poset has a maximal element) \cite{tukey40book}. Indeed, $A$ being a
non-empty inductive poset (see Lemma \ref{lem-a-mon} below), it has a
maximal element $a_k$. Since $a_k\le a_{k+1}$ and $a_k$ is maximal,
$a_{k+1}=a_k$. But the point of the previous authors was to prove
Zermelo's result that AC implies the well order theorem (every set can
be well ordered) \cite{zermelo04ma} without using ordinal theory (like
Zermelo in the second proof of his
theorem \cite{zermelo08ma}).\\

At the end of the 50's, by using Hartogs theorem (for any set $A$,
there is an ordinal $k$ that cannot be injected into $A$)
\cite{hartogs15ma}, Rubin and Rubin proved that, if $X$ is a strictly
inductive set of sets, $\le$ is inclusion and $f$ is extensive, then
$a_{k+1}=a_k$ for some $k$ \cite{rubin63book} (p. 18).\\

In 1957, Ward extended Tarski, Frink and Davis results to complete
semi-lattices (every non-empty subset has a $\lub$ but not necessarily
a $\glb$) \cite{ward57cjm}. Hence, a semi-lattice $X$ is complete iff,
for every $x\in X$, ${x\ad}=\{y\in X\mid y\le x\}$ is compact in the
interval topology; $\fp(f)$ is a complete semi-lattice if $X$ is a
complete semi-lattice and $f:X\a X$ is monotone; a semi-lattice is
compact in the interval topology iff every monotone $f:X\a X$ has a
fixed point.\\

In 1959, by refining Bourbaki's result, Abian and Brown extended
Tarski's result to strictly inductive posets having a pre-fixpoint of
$f$ \cite{abian61cjm}.\\

\comment{In 1960, Pelczar extended Tarski's result to any complete
  semi-lattice \cite{pelczar61apm}.}

In 1962, using Hartogs theorem, Devid\'e proved that, in a complete
lattice, if $f(x)=a_0\ou g(x)$ with $g$ monotone, then there is
$k\in\bO$ such that $a_{k+1}=a_k$ \cite{devide64fm} (note that $f$
does not need to be extensive, although it is so on $A$).\\

Up to now, we have seen that all conditions for $f$ to have a fixpoint
are requiring either $f$ to be monotone or $f$ to be extensive. So,
one may wonder what relations are known between these two classes of
functions, and whether one cannot devise a condition generalizing
both.

As for a relation between monotony and extensivity, we have:

\begin{lem}[\cite{bourbaki53tr189_nbr_092} p. 35]
\label{lem-mon-ext}
On a well ordered poset $(X,\le)$ (every non-empty subset has a least
element), any strictly monotone function $f:X\a X$ is extensive.
\end{lem}

\begin{prf}
Assume that $f$ is not extensive. Then, the set $E=\{x\in X\mid
{f(x)<x}\}$ is not empty. Let $e$ be its least element. By definition,
$f(e)<e$. By strict monotony, $f(f(e))<f(e)$. Hence, $f(e)\in E$ and,
$e\le f(e)$ by definition of $e$. Contradiction.\cqfd\\
\end{prf}

As for a condition subsuming both notions, we have:

In 1969, Salinas extended the previous fixpoint theorems
by requiring (P1) $a_0\in\pre(f)$, \ie $a_0\le f(a_0)$, and:

\begin{center}
(P2)\quad$f(x)\le f(y)$ if $x\le f(x)\le y$ \cite{salinas69tr}.
\end{center}

He provided two proofs, one using AC and Hartogs theorem, and another
one not using AC but some notion of chain\footnote{The notion of chain
  wrt a function $f$ has been first introduced by Dedekind for
  defining the notion of infinite set \cite{dedekind88book}. It has
  been used by Zermelo in his second proof that AC implies the well
  order theorem \cite{zermelo08ma}, and studied in details by
  Hessenberg \cite{hessenberg09crelle}.}. He also proved (using AC)
that, if $X$ is not strictly inductive but the set of upper bounds of
every non-empty chain of $X$ has a minimal element, $f$ satisfies (P1)
and (P2), and $\glb\{x,f(x)\}$ exists for all $x$, then $f$ has a
fixpoint.\\

In 1973, by using Bourbaki's theorem, Markowsky extended Tarski's
results to inductive (= chain-complete) posets (every chain has a
$\lub$), that is, $\fp(f)$ is chain-complete if $X$ is chain-complete,
and proved the converse, that is, $X$ is chain-complete if every
monotone (or $\glb$-preserving) map $f:X\a X$ has a {\em least}
fixpoint \cite{markowsky76au}.\\

In 1974, Pasini extended Salinas' result by proving that $\fp(f)$ has
a maximal element \cite{pasini74tr}.\\

In 1977, Cousot and Cousot studied the properties of
$a=(a_k)_{k\in\bO}$ when $f$ is monotone and $X$ a complete lattice
(remarking however that their results extend to posets every chain of
which has a $\lub$ and a $\glb$) \cite{cousot79pjm}. In particular,
they proved that, in a complete lattice, $A$ is bounded by every
post-fixpoint bigger than or equal to $a_0$, hence that $a$ can only
converge to the {\em least} fixpoint of $f$ bigger than or equal to
$a_0$, and indeed converges to this fixpoint if $f$ is monotone and
$a_0\le f(a_0)$. They also extended Tarski's result by showing that
$\pre(f)$ and $\post(f)$ are non-empty complete lattices too.

\section{Synthesis}

In conclusion, the most direct argument not using AC why there must be
some $k\in\bO$ such that $a_{k+1}=a_k$ is the one of Rubin
\cite{rubin63book} based on Hartogs theorem \cite{hartogs15ma}. We
hereafter split this result in two parts by showing first that, by
Hartogs theorem, $f$ has a fixpoint if $a$ is monotone, and second
that, $a$ is monotone if $f$ satisfies Salinas conditions (P2).

\begin{lem}
If $a$ is monotone, then there is an ordinal $k$ such that $a_{k+1}=a_k$.
\end{lem}

\begin{prf}
By Hartogs theorem, there is an ordinal $k$ that cannot be injected
into $A$ (the smallest such one is a cardinal). Therefore, $a|_k$, the
restriction of $a$ to $k$, is not an injection, that is, there are
$l_1<l_2<k$ such that $a_{l_1}=a_{l_2}$. Since $a$ is monotone, we
have $a_{l_1+1}=a_{l_1}$.\cqfd\\
\end{prf}

Now, one can easily prove that $a$ is monotone whenever $a_0$ is a
pre-fixpoint of $f$ and $f$ satisfies Salinas condition (P2) above:

\begin{lem}
\label{lem-a-mon}
$a$ is monotone if $a_0\le f(a_0)$ and $f(x)\le f(y)$ whenever $x\le
  f(x)\le y$, for all $x$ and $y$ in $A$.
\end{lem}

\begin{prf}
We prove that $a_k\le a_l$ whenever $k<l$ by induction on $l$ (1).
\begin{itemize}
\item If $l$ is a limit ordinal, then this is immediate.
\item Otherwise, $l=m+1$ and $k\le m$. If $k<m$ then, by induction
  hypothesis (1), we have $a_k\le a_m$. We now prove that, for all
  $i\le m$, $a_i\le a_{i+1}$, by induction on $i$ (2).
\begin{itemize}
\item $i=0$. $a_0\le a_1=f(a_0)$ by assumption.
\item $i=j+1$. By induction hypothesis (2), $a_j\le
  a_{j+1}=a_i$. Therefore, by (P2), $a_{j+1}=a_i\le a_{i+1}$.
\item $i=\lub(i)$. Let $j<i$. By induction hypothesis (2), $a_j\le
  a_{j+1}$. Since $j+1<i\le m$, by induction hypothesis (1),
  $a_{j+1}\le a_i$. Hence, by (P2), $a_{j+1}\le a_{i+1}$. Therefore,
  $a_j\le a_{i+1}$ and $a_i=\lub\{a_j\mid j<i\}\le a_{i+1}$.\cqfd
\end{itemize}
\end{itemize}
\end{prf}

Using a nice result of Abian and Brown \cite{abian61cjm} for {\em any}
poset $(X,\le)$ and {\em any} function $f$, we can go a little bit
further and instead consider the condition:

\begin{center}
(P2')\quad$f(x)\le f(y)$ if $x<f(x)\le y$ and there is no $z$ such
  that $x<z<f(x)$
\end{center}

\begin{dfn}[\cite{abian61cjm}]
A set $C\sle X$ is an {\em $a_0$-chain} if:
\begin{itemize}\itemsep=0mm
\item $C$ is well ordered;
\item $C$ has $a_0$ as least element;
\item $C$ is closed by non-empty $\lub$'s;
\item if $z\in C-\{\lub(C)\}$, then:
\begin{itemize}\itemsep=0mm
\item $f(z)\in C$,
\item $z<f(z)$,
\item there is no $y\in C$ such that $z<y<f(z)$.
\end{itemize}
\end{itemize}
Let $W$ be the set of elements $x\in X$ such that $x$ is the $\lub$ of
an $a_0$-chain.
\end{dfn}

In \cite{salinas69tr}, Salinas considered a similar (equal?) set
called the set of admissible subsets of $X$.

\begin{thm}[\cite{abian61cjm}]
\label{thm-abian}
For any poset $(X,\le)$, $a_0\in X$ and function $f:X\a X$:
\begin{itemize}\itemsep=0mm
\item $W$ is well ordered;
\item $W$ has $a_0$ as least element;
\item if $W$ has a $\lub$ $\xi$, then $W$ is an $a_0$-chain with $\xi$
  as greatest element and $\xi\not<f(\xi)$.
\end{itemize}
\end{thm}

Abian and Brown proved also that, for every $x\in W$, there is only
one $a_0$-chain $C$ such that $x=\lub(C)$, namely $\{y\in W\mid y\le
x\}$.

Note also that $W$ is not closed by $f$ in general. However, they
proved that, if $x\in W$ and $x\le f(x)$, then $f(x)\in W$.

\begin{thm}
\label{thm-fp}
In a non-empty strictly inductive poset $(X,\le)$, if $a_0\in\pre(f)$
and $f:X\a X$ is monotone on $W$, then $\lub(W)$ is a fixpoint of $f$.
\end{thm}

\begin{prf}
We simply follow the proof of Abian and Brown and check that, indeed,
the monotony of $f$ is used only on elements of $W$.

Since $X$ is strictly inductive, $\xi=\lub(W)$ exists. By the previous
theorem, $W$ is an $a_0$-chain and $\xi\not<f(\xi)$. Since $W$ has
$a_0$ as least element, we have $a_0\le\xi$. Since $\xi\not<f(\xi)$,
it suffices to check that $\xi\le f(\xi)$. If $a_0=\xi$, then this is
immediate since, by assumption, $a_0\le f(a_0)$. Assume now that
$a_0<\xi$. Then, since $W$ is an $a_0$-chain, we have $\xi\in W$ and
$V=W-\{\xi\}$ not empty, thus $\t=\lub(V)$ exists and
$\t\le\xi$. There are two cases:
\begin{itemize}
\item $\t=\xi$. Let $x\in V$. Then, $x<\xi$, $x<f(x)\in W$ and, by
  monotony of $f$ on $W$, $f(x)\le f(\xi)$. Hence, for all $x\in V$,
  $x<f(\xi)$. Therefore, $\xi\le f(\xi)$.
\item $\t<\xi=f(\t)$. Then, by monotony of $f$ on $W$, $\xi\le
  f(\xi)$.\cqfd\\
\end{itemize}
\end{prf}

Now, one can easily check that:

\begin{lem}
$f$ is monotone on $W$ if $f$ satisfies (P2').
\end{lem}

We now provide precise statements for Abian and Brown's claim that
${W=N}$. (Salinas also proved in \cite{salinas69tr} that his set of
{\em admissible} subsets of $X$ is $N$.)

\begin{lem}
$W\sle N$.
\end{lem}

\begin{prf}
If suffices to prove that $W$ is included in every set $Z$ containing
$a_0$ and closed by $f$ and non-empty $\lub$'s. We proceed by
well-founded induction on $<$. Let $x\in W$. If $x=a_0$, then we are
done. Assume now that $a_0<x$. After Lemma 4 in \cite{abian61cjm},
$x=\lub(C)$ with $C$ the $a_0$-chain $\{y\in W\mid y\le x\}$. We then
proceed as in Theorem \ref{thm-fp}. Since $C$ is an $a_0$-chain and
$a_0<x$, we have $D=C-\{x\}$ not empty, thus $\t=\lub(D)$ exists and
$\t\le x$. For all $y<x$, we have $y\in Z$ by induction
hypothesis. Therefore, $\t\in Z$ since $Z$ is closed by non-empty
$\lub$'s. If $\t=x$, then we are done. Otherwise, $x=f(\t)\in Z$ since
$Z$ is closed by $f$.\cqfd
\end{prf}

\begin{lem}
$N\sle W$ if $X$ is strictly inductive, $a_0\in\pre(f)$ and $f$ is
  monotone on $W$.
\end{lem}

\begin{prf}
It suffices to show that $W$ contains $a_0$ and is closed by $f$ and
non-empty $\lub$'s. By Theorem \ref{thm-abian}, we have $a_0\in W$ and
$W$ an $a_0$-chain. Hence, $W$ is closed by non-empty $\lub$'s. Let
$\xi=\lub(W)$ and $x\in W$. If $x=\xi$, then $f(x)=x\in W$ since
$\xi\in\fp(f)$ by Theorem \ref{thm-fp}. Otherwise, $x<\xi$ and
$f(x)\in W$ since $W$ is an $a_0$-chain.\cqfd\\
\end{prf}

For the sake of completeness, we also make precise Kuratowski's
relations between $A$ and $N$, when $X$ is strictly inductive (for $A$
to be well defined).

\begin{lem}[\cite{kuratowski22fm}]
$A\sle N$.
\end{lem}

\begin{prf}
It suffices to prove that $A$ is included in every set $Z\sle X$
containing $a$ and closed by $f$ and non-empty $\lub$'s, by
transfinite induction. Let $a_k\in A$. If $k=0$, then $a_k\in Z$ by
assumption. If $k=j+1$, then $a_k=f(a_j)\in Z$ since, by induction
hypothesis, $a_j\in Z$ and $Z$ is closed by $f$. Finally, if
$k=\lub(k)$, then $a_k=\lub\{a_j\mid j<k\}\in Z$ since $Z$ is closed
by non-empty $\lub$'s and, for all $j<k$, $a_j\in Z$ by induction
hypothesis.\cqfd
\end{prf}

\begin{lem}
$N\sle A$ if $a$ is monotone.
\end{lem}

\begin{prf}
Since $A$ contains $a_0$ and is closed by $f$, it suffices to prove
that $A$ is closed by non-empty $\lub$'s. Let $Z$ be a non-empty
subset of $A$. Then, there is a set $K$ of ordinals such that
$Z=\{x\mid\ex k\in K,x=a_k\}$. Since $a$ is monotone, $\lub(Z)$ exists
and equals $a_k$ where $k=\lub(K)$ (every {\em set} of ordinals has a
$\lub$).\cqfd\\
\end{prf}

Hence, we can conclude:

\begin{thm}
If $X$ is strictly inductive, $a_0\in\pre(f)$ and $f$ satisfies (P2'),
then $N=W=A$. Therefore, $N$, $W$ and $A$ are $a_0$-chains and there
is $k$ such that $a_{k+1}=a_k=\lub(N)$ is the {\em least} fixpoint of
$f$ bigger than or equal to $a_0$.
\end{thm}

\begin{prf}
Since $f$ satisfies (P2'), $f$ is monotone on $W$. Hence, $N=W$. Since
$A\sle N$, $f$ is monotone on $A$. Hence, $a$ is monotone. Therefore,
$N=A$.\cqfd\\
\end{prf}

Note that, if $\xi=\lub(W)$, then $f$ is both strictly monotone and
strictly extensive on $W-\{\xi\}$.

We finish with some ultimate remarks.

In \cite{devide64fm}, with $X$ a complete lattice, Devid\'e takes
$f(x)=a_0\ou g(x)$ with $g$ monotone. In this case, one can easily
check that, if $a_0\le g(a_0)$, then $f$ and $g$ have the same set of
fixpoints bigger than or equal to $a_0$. Moreover, the transfinite
iterates of $f$ and $g$ are equal.

Now, consider $f(x)=x\ou g(x)$ with $g$ monotone. Then, $f$ is both
monotone and extensive, and $f$ and $g$ are equal on
$\pre(g)$. Moreover, $\fp(g)\sle\fp(f)=\post(g)$ and, if $X$ has at
least two distinct elements, then $\fp(f)\not\sle\fp(g)$: if $g$ is
the constant function equal to the least element $\bot$ of $X$, then
$\fp(g)=\{\bot\}$ and $\fp(f)=X\neq\{\bot\}$ since $f$ is the identity
and $X$ has at least two elements. However, the {\em least} fixpoint
of $f$ is also the {\em least} fixpoint of $g$. Moreover, if $a_0\le
g(a_0)$, then the transfinite iterates of $f$ and $g$ are equal.

\bibliographystyle{alpha}

\input{biblio}
\end{document}

%% file: macros.tex
\newcommand{\comment}[1]{}

\newcommand{\vide}{\emptyset}

\newcommand{\ie}{{\em i.e.} }

\newcommand{\glb}{\mr{glb}}
\newcommand{\lub}{\mr{lub}}

\renewcommand{\a}{\rightarrow}

\newcommand{\ad}{\downarrow}

\newcommand{\ex}{\exists}

\newcommand{\ou}{\vee}

\newcommand{\sle}{\subseteq}

\renewcommand{\prod}{_\mr{prod}}

\renewcommand{\t}{\theta}

\newcommand{\w}{\omega}

\newcommand{\mr}{\mathrm}
\newcommand{\mb}{\mathbb} 

\newcommand{\bO}{\mb{O}}

\newcommand{\fp}{\ms{p}}

\newenvironment{rul}
  {$\begin{array}{rcl}}
  {\end{array}$}

\newenvironment{rew}[1][~~\a~~]
  {$\begin{array}{r@{#1}l}}
  {\end{array}$}
\newenvironment{rewc}[1][~~\a~~]
  {\begin{center}\begin{rew}[#1]}
  {\end{rew}\end{center}}

\newcounter{dfnnum}
\newcounter{lemnum}
\newcounter{thmnum}

{\theorembodyfont{\rmfamily} 
  \newtheorem{dfn}[dfnnum]{Definition}
  \newtheorem{lem}[lemnum]{Lemma}
  \newtheorem{thm}[thmnum]{Theorem}

}

\newcommand{\cqfd}{\hfill$\blacksquare$}

\newenvironment{prf}{{\bf Proof.}}{}

\newenvironment{lstgeneric}[2]
  {\begin{list}{#1}{\topsep=.5mm\itemsep=.5mm\parsep=0mm%
    \itemindent=-3ex\labelsep=1ex\labelwidth=0ex #2}}
  {\end{list}}
